\documentclass[11pt]{article}

 \usepackage{amsmath}
 \usepackage{amssymb}
 \usepackage{graphicx}
 \usepackage{verbatim}
  
 \newtheorem{thm}{Theorem}[section]

 \newtheorem{prop}[thm]{Proposition}

 \newtheorem{lemma}[thm]{Lemma}
 \newtheorem{rem}[thm]{Remark}

 \title{Convexity of the median in the gamma distribution}
 \author{Christian Berg and Henrik L. Pedersen\footnote{Research supported by the
     Carlsberg Foundation}}

 \date{\today}

\begin{document}

 \maketitle

 \begin{abstract} We show that the median $m(x)$ in the gamma distribution
   with parameter $x$ is a strictly convex function on the positive half-line.
 \end{abstract}

\noindent 
2000 {\em Mathematics Subject Classification}:\\
primary 60E05; secondary 41A60, 33B15. 

\noindent
Keywords: median, gamma function, gamma distribution.
\bigskip

\section{Introduction}
The median of the gamma distribution with (positive) parameter $x$ is
defined implicitly by the formula 
\begin{equation}
\label{eq:def_m}
\int_0^{m(x)}e^{-t}t^{x-1}\, dt=\frac{1}{2} \int_0^{\infty}e^{-t}t^{x-1}\, dt.
\end{equation}
In a recent paper (see \cite{BP}) we showed the $0<m'(x)<1$ for all
$x>0$. Consequently, $m(x)-x$ is a decreasing function, which for
$x=1,2,\ldots$ yields a positive answer to the Chen-Rubin
conjecture. Other authors have solved this conjecture in its discrete
setting (see \cite{Alm}, \cite{Adell}, \cite{A}).

In \cite{A2} convexity of the sequence $m(n+1)$ has been established, and
the natural question arises if $m(x)$ is a convex function. The main
result of this paper is the following.
\begin{thm}
\label{thm:main}
The median $m(x)$ defined in (\ref{eq:def_m}) satisfies $m''(x)>0$. In
particular it is a strictly convex function for
$x>0$.
\end{thm}

\section{Proofs}
The proof is based on some results in \cite{BP}, which we briefly
describe.
Convexity of $m$ is studied through the function 
\begin{equation}\label{eq:phi}
\varphi(x)\equiv \log\frac{x}{m(x)}, \quad{x>0}.
\end{equation}
This function played a key role in \cite{BP}, and we recall its crucial
properties in the proposition below.
\begin{prop}
\label{prop:xphi}
The function $x\to x\varphi(x)$ is strictly decreasing for $x>0$ and 
\begin{eqnarray*}
\lim_{x\to 0_+}x\varphi (x) & = & \log 2,\\
\lim_{x\to \infty}x\varphi (x) & = & \frac{1}{3}.
\end{eqnarray*}
\end{prop}
\begin{rem} Proposition \ref{prop:xphi} is established by showing
 $(x\varphi(x))'<0$. It follows that the function
  $\varphi(x)$ is itself strictly decreasing and
  $\varphi(x)<-x\varphi'(x)$.
\end{rem}
The starting point for proving Theorem \ref{thm:main} is the relation
\cite[(10)]{BP}
\begin{equation}
\label{eq:gAB}
(x\varphi(x))'=-e^{g(x)}(A(x)+B(x)),
\end{equation}
where
\begin{eqnarray*}
g(x)&=&x(\varphi(x)-1+e^{-\varphi(x)}),\\
A(x)&=& \int_{0}^{x\varphi (x)}e^{-s}e^{x(1-e^{-s/x})}\left(
  1-\left( 1+\frac{s}{x}\right) e^{-s/x}\right)\, ds,\\
B(x)&=& \frac{1}{2}\int_{0}^{\infty}te^{-xt}\xi'(t+1)\, dt,
\end{eqnarray*}
and where $\xi$ is a certain positive, increasing and concave function
on $[1,\infty)$ satisfying $\xi'(t+1)<8/135$ for $t>0$. To establish
these properties of $\xi$ is quite involved, and we refer to
\cite[Section 5]{BP} for details.

Before proving the theorem we state the following lemmas, whose proofs
are given later.
\begin{lemma}
\label{lemma:g}
For the function $g$ we have 
\begin{eqnarray*}
g(x) & < & x\varphi(x),\\
-g'(x) & <& -x\varphi'(x)\varphi(x) \quad \text{and}\\
-g'(x) & <& -x\varphi'(x) 
\end{eqnarray*}
for all $x>0$.
\end{lemma}

\begin{lemma}
\label{lemma:A}
We have for $x>0$
\begin{eqnarray*}
A(x) & < & \frac{x\varphi(x)^3}{6},\\
-A'(x) & < & -\frac{1}{6}\varphi(x)^3-\frac{1}{2}x\varphi'(x)\varphi(x)^2.
\end{eqnarray*}
\end{lemma}

\begin{lemma}
\label{lemma:B}
We have for $x>0$
\begin{eqnarray*}
B(x) & < & \frac{4}{135x^2},\\
-B'(x) & < & \frac{8}{135x^3}.
\end{eqnarray*}
\end{lemma}

{\it Proof of Theorem \ref{thm:main}.} From equation (\ref{eq:phi}) we get
$$
m''(x)=-e^{-\varphi(x)}\left(2\varphi'(x)+x\varphi''(x)-x\varphi'(x)^2\right),
$$
so that $m''(x)>0$ is equivalent to the inequality
$$
(x\varphi(x))''<x\varphi'(x)^2.
$$
Differentiation of (\ref{eq:gAB}) yields
$$
(x\varphi(x))''=e^{g(x)}(-g'(x))(A(x)+B(x))+e^{g(x)}(-A'(x)-B'(x)).
$$
By using Lemma \ref{lemma:A} and \ref{lemma:B} it follows that
\begin{eqnarray*}
-A'(x)-B'(x) &< &
-\frac{1}{2}x\varphi'(x)\varphi(x)^2+\frac{8}{135x^3}-\frac{1}{6}\varphi(x)^3\\
&= &-\frac{1}{2}x\varphi'(x)\varphi(x)^2+\frac{\varphi(x)^2}{x}\left(\frac{8}{135(x\varphi(x))^2}-\frac{1}{6}x\varphi(x)\right).
\end{eqnarray*}
Here the expression
in the brackets is positive,
since $(x\varphi(x))^3< (\log 2)^3< 48/135$. Therefore,
and because $\varphi(x)<-x\varphi'(x)$,
\begin{eqnarray*}
-A'(x)-B'(x) &< &\frac{1}{2}x\varphi'(x)^2\,
x\varphi(x)+x\varphi'(x)^2\left(\frac{8}{135(x\varphi(x))^2}-\frac{1}{6}x\varphi(x)\right)\\
& = & x\varphi'(x)^2 \left(\frac{8}{135(x\varphi(x))^2}+\frac{1}{3}x\varphi(x)\right).
\end{eqnarray*}
We also have from Lemma \ref{lemma:g}, \ref{lemma:A} and \ref{lemma:B}, 
\begin{eqnarray*}
-g'(x)(A(x)+B(x)) &< &
-x\varphi'(x)\varphi(x)\left(\frac{x\varphi(x)^3}{6}+\frac{4}{135x^2}\right)
\\
& < &
x^2\varphi'(x)^2\varphi(x)^2\left(\frac{x\varphi(x)}{6}+\frac{4}{135(x\varphi(x))^2}\right).
\end{eqnarray*}
Combination of these inequalities yields
\begin{eqnarray*}
(x\varphi(x))'' &< &
x\varphi'(x)^2e^{x\varphi(x)}\left(x\varphi(x)^2\left(\frac{x\varphi(x)}{6}+\frac{4}{135(x\varphi(x))^2}\right)\right.
    \\
&& \left. +\frac{8}{135(x\varphi(x))^2}+\frac{1}{3}x\varphi(x)\right).
\end{eqnarray*}

Supposing that $x\geq 1$, it follows that 
\begin{eqnarray*}
(x\varphi(x))'' &< &
x\varphi'(x)^2e^{x\varphi(x)}\left((x\varphi(x))^2\left(\frac{x\varphi(x)}{6}+\frac{4}{135(x\varphi(x))^2}\right)\right.\\
&&\left. +\frac{8}{135(x\varphi(x))^2}+\frac{1}{3}x\varphi(x)\right)\\
&=&x\varphi'(x)^2e^{x\varphi(x)}\left(\frac{(x\varphi(x))^3}{6}+\frac{4}{135}+\frac{8}{135(x\varphi(x))^2}+\frac{1}{3}x\varphi(x)\right)\\
&=&x\varphi'(x)^2 h_1(x\varphi(x)),
\end{eqnarray*}
where $h_1$ is given by
$$
h_1(t)= e^{t}\left(\frac{t^3}{6}+\frac{4}{135}+\frac{8}{135t^2}+\frac{t}{3}\right).
$$
One can show that $h_1$ attains its maximum on the interval $[1/3,
\log 2]$ at the left end point and that
$h_1(1/3)=(551/810)\sqrt[3]{e}\approx 0.9494$. Therefore it follows that $(x\varphi(x))'' <x\varphi'(x)^2$ for
$x\geq 1$.

For $0<x<1$ the estimate $-g'(x)<-x\varphi'(x)$ from Lemma \ref{lemma:g}
is used and in this way
we get
$$
(x\varphi(x))'' <x\varphi'(x)^2 h_2(x\varphi(x)),
$$
where 
$$
h_2(t)= e^{t}\left(\frac{t^2}{6}+\frac{4}{135t}+\frac{8}{135t^2}+\frac{t}{3}\right).
$$
Since $x<1$ and $x\varphi(x)$ decreases we have $x\varphi(x)>
\varphi(1)=-\log \log 2$. One can show that $h_2$ attains its maximum on
the interval $[-\log \log 2, \log 2]$ for $t=-\log \log 2$ and
that $h_2(-\log \log 2)\approx 0.9616$. Therefore $(x\varphi(x))''
<x\varphi'(x)^2$ for $x< 1$.  \hfill $\square$

\begin{rem}
  The function $h_2$ becomes larger than 1 on the interval $[1/3,\log
  2]$, so $h_2$ cannot be used to obtain the
  inequality $(x\varphi(x))'' <x\varphi'(x)^2$ for all $x>0$.
\end{rem}

{\it Proof of Lemma \ref{lemma:g}.} It is clear that
$g(x)<x\varphi(x)$. Differentiation yields
\begin{eqnarray*}
-g'(x) & = & -\varphi(x)+(1-x\varphi'(x))(1-e^{-\varphi(x)})\\
&<&-\varphi(x)+(1-x\varphi'(x))\varphi(x)\, = \,
-x\varphi'(x)\varphi(x),
\end{eqnarray*}
where we have used  $1-e^{-a}<a$ for $a>0$.

To find an estimate that is more accurate for $x$ near 0 we use
\begin{eqnarray*}
-g'(x) & = & -x\varphi'(x)(1-e^{-\varphi(x)})-\varphi(x)+1-e^{-\varphi(x)}\\
&<&-x\varphi'(x)-\varphi(x)+1-e^{-\varphi(x)}\, < \,-x\varphi'(x).
\end{eqnarray*}
\hfill $\square$

{\it Proof of Lemma \ref{lemma:A}.}
Using that $1-e^{-a}<a$ and $1-(1+a)e^{-a}< a^2/2$  for
$a> 0$, we can estimate $A(x)$ by 
$$
A(x)< \int_{0}^{x\varphi (x)}e^{-s}e^{x(s/x)}\left(
  \frac{s^2}{2x^2}\right)\, ds =
\frac{x\varphi(x)^3}{6}.
$$
A computation shows that
\begin{eqnarray*}
-A'(x)&=&
-(\varphi(x)+x\varphi'(x))e^{-x\varphi(x)}e^{x(1-e^{-\varphi(x)})}(1-(1+\varphi(x))e^{-\varphi(x)})\\
&&-\int_{0}^{x\varphi (x)}e^{-s}e^{x(1-e^{-s/x})}\left(
  1-\left( 1+\frac{s}{x}\right) e^{-s/x}\right)^2\, ds\\
&&+\int_{0}^{x\varphi
  (x)}e^{-s}e^{x(1-e^{-s/x})}\frac{s^2}{x^3}e^{-s/x}\, ds\\
&<&-(\varphi(x)+x\varphi'(x))\frac{1}{2}\varphi(x)^2+\int_{0}^{x\varphi
  (x)}s^2e^{-s/x}\, ds\frac{1}{x^3}\\
&<&
-\frac{1}{2}\varphi(x)^3-\frac{1}{2}x\varphi'(x)\varphi(x)^2+\frac{1}{3}
\varphi(x)^3\\
&=&-\frac{1}{6}\varphi(x)^3-\frac{1}{2}x\varphi'(x)\varphi(x)^2.
\end{eqnarray*}
\hfill $\square$

{\it Proof of Lemma \ref{lemma:B}.} These estimates follow directly
from the inequality $\xi'(t+1)< 8/135$.\hfill $\square$

\author{Department of Mathematics, University of Copenhagen,
  Universitetsparken 5, DK-2100, Copenhagen, Denmark.}\\
\noindent
\author{Email: berg@math.ku.dk\\Fax: +4535320704}
 \bigskip

\noindent
\author{Department of Natural Sciences, Royal Veterinary and
  Agricultural University,
  Thorvaldsensvej 40, DK-1871, Copenhagen, Denmark.}\\
\noindent
\author{Email: henrikp@dina.kvl.dk\\Fax: +4535282350}

\end{document}